\documentstyle[11pt]{article}
\title{{\bf Comment on a Paper by Yucai Su On Jacobian Conjecture (Dec 30, 2005)}}
\author{{\Large T.T. Moh}\thanks{ Math Department, Purdue University, West Lafayette, Indiana 47907-1395. tel: (765)-494-1930, e-mail ttm@math.purdue.edu}}
\begin{document}
\maketitle

\begin{abstract}
The said paper [Su2] entitled "Proof Of Two Dimensional Jacobian Conjecture" is false.
\end{abstract}

\noindent{\bf Comments}

The Jacobian Conjecture is a well-known hard problem. Recently there are many attempts
to solve it. This paper is one of the false attempts.

\noindent{\bf Minor mistakes}

We have read Su's latest version (Dec 30, 2005) of his papers. There are numerous typos (say, some $d's$,
and all $d_i's$ really should be $d_2$ etc.) and unconventional concepts and notations (say, "rational 
polynomials (it is a CS notation, not commonly used in Mathematics)"(=rational functions?, or =polynomials 
with rational numbers as coefficients?) "modulo {\bf F}(x,y)" means
"{\bf F}[x](($y^{-1}$)) as a ring or a commutative group?". Note that in those two cases, the
meaning of "modulo {\bf F}(x,y)" is totally different.) Those errors should be fixed. 

\noindent{\bf Major mistakes}

One of the mistakes is that on pg 17 of [Su2], let us quote, 

{\it Lemma 3.10: a'=1. In particular p'=1, q'=$d_2$ namely, $\bar{F}_1$(x,y)=y$^{d_2}+c_{d_2}x$ for some
$c_{d_2}\in F^{*}$}.

There are two problems about this lemma which is the kernel of the whole proof in his paper (all previous
statements in this paper are largely known in different terminologies) , (1) the proof 
of the lemma is with gaps, (2) the lemma can not be true since there are counter-examples.

(1). The said lemma is based on formula (3.56) on page 16. Using formula (3.56) he makes several computations,
and he uses the argument (on page 17, the 8th line from the bottom) that {\it But the sum of these two terms 
is nonzero (with coefficient $\alpha_1i_k+\alpha_2j_k<0)$}. From his previous arguments, we only know that
$$
\alpha"i_k+1-j_k>0
$$
which can be see as

\begin{eqnarray*}
&&\alpha"i_k+1-j_k>0\Leftrightarrow 1>-\alpha"i_k+j_k\\
&&\Leftrightarrow d_2p'>-(\alpha"d_2p')i_k+(d_2p')j_k\\ 
&&\Leftrightarrow d_2p'-(p'+q')i_k>-(\alpha"d_2p'+p'+q')i_k+(d_2p')j_k\\
&&\Leftrightarrow d_2p'-(p'+q')i_k>(\alpha_1)i_k+(d_2p')j_k
\end{eqnarray*}
There is no way to show that
$$
0>d_2p'-(p'+q')i_k
$$
Therefore he can not use $\alpha_1i_k+\alpha_2j_k<0$. Hence his conclusion is false.

(2). BTW, as we indicated in [M2] that {\it "For the last thirty years, the only
case one can not handle is the case of more than one point at $\infty$, i.e., in Su's notation,
$p=1,j=m+n-2$ (cf pg 17, Subcase of [Su1]), all other cases are well-known. The author should give a
convincing argument for this case (in fact, for this case only)"}. Namely he should only prove his
"{\it Subcase 2.2.1 (on page 16 of [Su2]): Suppose $B\not=\emptyset$ and j=min B satisfies (m+n-1-j)p=1
(this is the most nontrivial case)"} with p=1. 

We may find an example to formula (3.56) with $p'=1,q'=1, (p=p'/q'=1), d_2=3, \alpha'=0, \bar{F}_1=y^2(y+x), 
P=y(y+x),c_2'=-x$.

Note that the formula (3.56) is reduced to

$$
-2\bar{F}_1^{-1}P\partial_y \bar{F}_1+3\partial_y P= -x
$$
and satisfied by our data.

Thus the conclusion $\alpha'=1$ of lemma 3.10 is false.

\noindent{\bf A few words for Mr. Su}

The problem of Jacobian Conjecture is very hard. Perhaps it will take human being another 100 years to
solve it. Your attempt is noble, Maybe the Gods of Olympus will smile on you one day. Do not be
too disappointed. B. Sagre has the honor of publishing three wrong proofs and C. Chevalley mistakes a
wrong proof for a correct one in the 1950's in his Math Review comments, and I.R. Shafarevich uses Jacobian
Conjecture (to him it is a theorem) as a fact. You are in a good company. One only remembers the correct
statements from Scientists and Mathematicians, nobody remembers the wrong ones. 
 
\noindent{\bf Reference} 

[M1] {\bf Moh, T.T.}:  {\it On the Jacobian Conjecture} Crelle 1983.

[M2] {\bf Moh, T.T.}:  {\it  Comment on a Paper by Yucai Su On Jacobian Conjecture (Dec 18, 2005) Arxiv,
Dec 20, 2005. 

[Su1] {\bf Su, Yucai}  {\it Proof Of Two Dimensional Jacobian Conjecture} Arxiv, Dec 18, 2005,

[Su2] {\bf Su, Yucai}  {\it Proof Of Two Dimensional Jacobian Conjecture} Arxiv, Dec 30, 2005,

\end{document}